\documentclass[11pt]{amsart}

\usepackage{amsmath,amsfonts,amssymb,amscd}
\begin{document}
\renewcommand{\thefootnote}{\fnsymbol{footnote}}
\pagestyle{plain}

\title{The Yau-Tian-Donaldson Conjecture  \\
for general polarizations}
\author{Toshiki Mabuchi${}^*$}
\maketitle
\footnotetext{ ${}^{*}$Supported 
by JSPS Grant-in-Aid for Scientific Research (A) No. 20244005.}
\abstract
In this paper, assuming that a polarized algebraic manifold $(X,L)$ is strongly K-stable  in the sense of \cite{M},
we shall show that the 
class $c_1(L)_{\Bbb R}$ admits a constant scalar curvature K\"ahler metric.
Since strong K-stability 
implies asymptotic Chow-stability (cf.~\cite{MN1}), 
we have a sequence $\{\omega_i\}$ of balanced metrics in the class $c_1(L)_{\Bbb R}$.
Replace the sequence by its suitable subsequence if necessary.
Then if $\{\omega_i\}$ were not convergent, the associated sequence $\{\mu_i\}$ of polarized test configurations would satisfy the inequality
$$
F_1 (\{\mu_i\}) \geq 0 
$$ 
in contradiction to strong K-stability for $(X,L)$. Hence the sequence $\{\omega_i\}$ converges to a constant scalar curvature K\"ahler metric in $c_1(L)_{\Bbb R}$.
\endabstract

\section{Introduction}

By a {\it polarized algebraic manifold} $(X,L)$, we mean
a pair of a nonsingular irreducible projective algebraic variety $X$, 
defined over $\Bbb C$, and a very ample line bundle $L$
over $X$. 
Replacing $L$ by its positive integral multiple if necessary, we may assume that
$$
H^q(X, L^{\otimes\ell}) = \{0\}, \qquad \ell = 1,2,\dots ;\,
q= 1,2,\dots,n,
$$
where $n$ is the complex dimension of $X$.
In this paper,
we fix once for all such a pair $(X,L)$. 
For the affine line 
 $ \Bbb A^1 := \{z\in \Bbb C\}$,
let the algebraic torus $T := \Bbb C^*$ act on $\Bbb A^1$ by multiplication of complex numbers
$$
T\times \Bbb A^1 \to \Bbb A^1,
\qquad (t, z) \mapsto
t z.
$$
By fixing a Hermitian metric $h$ for $L$ such that $\omega := c_1(L;h)$ is K\"ahler, we endow the space $V_{\ell}:= H^0(X,L^{\otimes\ell})$ 
of holomorphic sections for $L^{\otimes \ell}$ with the Hermitian metric $\rho_{\ell}$ defined by
$$
\langle\sigma', \sigma'' \rangle_{\rho^{}_{\ell}} := \int_X \; (\sigma',\sigma'' )^{}_h \; \omega^n,
\qquad \sigma',\sigma'' \in V_{\ell},
$$
where $(\sigma',\sigma'' )_h$ denotes the pointwise Hermitian inner product of 
$\sigma'$ and $\sigma''$ by the $\ell$-multiple of $h$.
For the Kodaira embedding $
 \Phi_{\ell} \,:\, X \, \hookrightarrow \, \Bbb P^*(V_{\ell})$
associated to the complete linear system $|L^{\otimes\ell}|$ on $X$, 
we put $X_{\ell} := \Phi_{\ell}(X)$.
Let
$\psi \, : \,\Bbb C^* \,\to \,\operatorname{GL}(V_{\ell})$
be an algebraic group homomorphism
such that the compact subgroup $S^1 \subset \Bbb C^* \,(= T)$ acts isometrically on $(V_{\ell}, \rho_{\ell})$.
Take the irreducible algebraic subvariety $\mathcal{X}^{\psi}$ of $\Bbb A^1 \times \Bbb P^* (V_{\ell})$ obtained as 
the closure of $\cup_{z\in \Bbb C^*} \mathcal{X}^{\psi}_z$ in $\Bbb A^1 \times \Bbb P^* (V_{\ell})$. 
Here we set
$$
\mathcal{X}^{\psi}_z := \{z\}\times\psi (z)  \Phi_{\ell} (X),
\qquad z \in \Bbb C^*,
$$
and $\psi (z)  $ in $\operatorname{GL}(V_{\ell})$ acts naturally on the space  
$\Bbb P^* (V_{\ell})$ of all hyperplanes in $V_{\ell}$ passing through the origin.
We then consider the map
$$
\pi : \mathcal{X}^{\psi} \to \Bbb A^1
$$ 
induced by the projection of $\Bbb A^1 \times \Bbb P^* (V_{\ell})$ to the first factor $\Bbb A^1$.
Moreover, for the hyperplane bundle $\mathcal{O}_{\Bbb P^*(V_{\ell})}(1)$ on $\Bbb P^*(V_{\ell})$, 
we consider the pullback 
$$
\mathcal{L}^{\psi}\, :=\,\operatorname{pr}_2^*\mathcal{O}_{\Bbb P^*(V_{\ell})}(1)_{|\mathcal{X}^{\psi}},
$$
where $\operatorname{pr}_2 : \Bbb A^1 \times \Bbb P^* (V_{\ell}) \to \Bbb P^* (V_{\ell})$
denotes the projection to the second factor.
For the dual space $V_{\ell}^*$ of $V_{\ell}$,
the $\Bbb C^*$-action on $\Bbb A^1 \times V_{\ell}^*$ defined by
$$
\Bbb C^* \times (\Bbb A^1 \times V_{\ell}^*)\to \Bbb A^1 \times V_{\ell}^*,
\quad (t, (z, p))\mapsto  (tz, \psi (t) p),
$$
naturally induces the $\Bbb C^*$-action on $\Bbb A^1 \times \Bbb P^*(V_{\ell})$ and $\mathcal{O}_{\Bbb P^*(V_{\ell})}(-1)$, where $\operatorname{GL}(V_{\ell})$ acts on $V_{\ell}^*$ by the contragradient representation.
 This then induces $\Bbb C^*$-actions on $\mathcal{X}^{\psi}$ and $\mathcal{L}^{\psi}$, and $\pi : \mathcal{X}^{\psi} \to \Bbb A^1$ is a $\Bbb C^*$-equivariant 
 projective morphism with relative very ample line bundle 
$\mathcal{L}^{\psi}$ such that
$$
(\mathcal{X}^{\psi}_z, \mathcal{L}_z^{\psi})\; \cong \; (X,L^{\otimes \ell}),
\qquad z \neq 0,
$$
where $\mathcal{L}_z^{\psi}$ is the restriction of $\mathcal{L}^{\psi}$ to $\mathcal{X}^{\psi}_z := \pi^{-1}(z)$.
Then a  triple $({\mathcal{X}}, {\mathcal{L}}, \psi )$ is called 
a {\it test configuration for $(X,L)$}, if we have both $\mathcal{X} = \mathcal{X}^{\psi}$ and  $\mathcal{L}= \mathcal{L}^{\psi}$.
Here $\ell$ is called the {\it exponent} of $({\mathcal{X}}, {\mathcal{L}}, \psi )$.
From now on until the end of  Step 1 of Section 4, for 
$(\mathcal{X},\mathcal{L},\psi )$ to be a test configuration,
we make an additional assumption that
$\psi$  is  written
in the form
$$
\psi : \Bbb C^* \to  \operatorname{SL} (V_{\ell}).
$$
Then $({\mathcal{X}}, {\mathcal{L}}, \psi )$ is called {\it trivial},
if $\psi$ 
is a trivial homomorphism.
We now consider the set $\mathcal{M}$ of all sequences $\{\mu_j \}$ of test configurations 
$$
\mu_j  \, =\, (\mathcal{X}_j, \mathcal{L}_j, \psi_j ), \qquad j =1,2,\dots,
$$
for $(X,L)$  
such that for each $j$, the exponent $\ell_j$ of the test configuration $\mu_j$ satisfies the following 
condition:
$$
\text{$\ell_j \to \infty$,  \;\; as $j \to \infty$.}
$$
In \cite{M}, for each $\{\mu_j\}\in  
\mathcal{M}$, we defined the Donaldson-Futaki invariant 
$F_1 (\{\mu_j\}) \in  \Bbb R \cup \{-\infty\}$.
Then we  have the strong version of 
K-stability and K-semistability as follows:

\medskip\noindent
{\em Definition\/ $1.1$}. 
(1) The polarized algebraic manifold $(X,L)$ is called 
{\it strongly $K$-semistable}, 
if $F_1(\{\mu_j\} ) \leq 0$ for all 
$\{\mu_j\} \in 
\mathcal{M}$.

\medskip\noindent
(2) A strongly K-semistable polarized algebraic manifold $(X,L)$
 is called 
{\it strongly $K$-stable},
if  for every $\{\mu_j\}\in\mathcal{M}$ satisfying $ F_1(\{\mu_j\} ) = 0$,
there exists a $j_0$ such that
 $\mu_j$ are trivial for all $j$ with $j \geq j_0$.
 
\medskip
Recall that these stabilities are independent of the choice 
of the Hermitian metric $h$ for $L$ (see \cite{MN2}).
The purpose of this paper is to show the following:

\medskip\noindent
{\bf Main Theorem.} 
{\em  If $(X,L)$ is strongly K-stable, then
the class $c_1(L)_{\Bbb R}$ admits a constant scalar curvature K\"ahler metric.}

\section{The Donaldson-Futaki invariant $F_1$ on $\mathcal{M}$}

\medskip\noindent
{\em Definition \/$2.1$.} For a complex vector space $V$,
let $\phi : T \to \operatorname{GL}(V)$ be an algebraic group homomorphism. For the real Lie subgroup 
$$
T_{\Bbb R}\; :=\;\{t \in T\,;\, t\in \Bbb R_+\}
$$ 
of the algebraic torus $T = \{t \in \Bbb C^*\}$,
we define the associated Lie group homomorphism 
$\phi^{\operatorname{SL}} : T_{\Bbb R} \to \operatorname{SL}(V)$ by
$$
 \phi^{\operatorname{SL}} (t) \;:=\; \frac{\phi (t) }{\det (\phi (t))^{1/N}},
 \qquad t \in T_{\Bbb R},
 $$
where $N:=\dim V$. Let $b_1$, $b_2$, \dots , $b_N$ be the weights of the action
by $ \phi^{\operatorname{SL}}$ on the dual vector space $V^*$ of $V$, so that
we have the equalities 
$$
\phi^{\operatorname{SL}} (t) \cdot \sigma_{\alpha} \, =\, t^{-b_{\alpha}}\sigma_{\alpha},
\qquad \alpha = 1,2,\dots,N,
$$
for some basis $\{\sigma_1, \sigma_2, \dots, \sigma_N\}$ 
of $V$.
Then we define $\|\phi\|_1$ and $\|\phi\|_{\infty}$ by
$$
\|\phi \|_1 \; := \Sigma_{\alpha = 1}^N \, |b_{\alpha} | 
\quad \text{ and }\quad
\|\phi\|_{\infty}:= \max \{\, |b_1|, |b_2|, \dots, |b_N|\,\}.
$$

\medskip\noindent
{\em Definition \/$2.2$.}
Put $d := \ell^{\,n} c_1(L)^n [X]$. For $(V_{\ell}, \rho_{\ell})$ in the introduction, we define
a  space $W_{\ell}$ by
$$
W_{\ell} \, := \;\{\operatorname{Sym}^{d}(V_{\ell})\}^{\otimes n+1},
$$
where  $\operatorname{Sym}^{d}(V_{\ell} )$ is the $d$-th symmetric tensor product  of $V_{\ell}$.
Then the dual space $W_{\ell}^*$ of $W^{}_{\ell}$ admits the Chow norm (cf.~\cite{Zh})
$$
W_{\ell}^* \owns w \;\mapsto \;\| w\|^{}_{\operatorname{CH}(\rho^{}_{\ell})} \in \Bbb R^{}_{\geq 0},
$$
associated to the Hermitian metric $\rho^{}_{\ell}$ on $V_{\ell}$. For the Kodaira embedding
$\Phi_{\ell}: X \hookrightarrow  \Bbb P^*(V_{\ell})$ as in the introduction, let
$$
0 \neq \hat{X}_{\ell}\in W_{\ell}^*
$$ 
be the associated Chow form for $X_{\ell} = \Phi_{\ell}(X)$ 
viewed as an irreducible reduced algebraic cycle 
 on the projective space $\Bbb P^*(V_{\ell})$.

\bigskip
Let $\mu_j = (\mathcal{X}_j, \mathcal{L}_j, \psi_j)$, $ j = 1,2,\dots$,
be a sequence of test configurations 
for $(X,L)$. 
We then define $\|\mu_j\|_1$ and $\|\mu_j\|_{\infty}$ by
$$
 \|\mu_j\|_1 \; := \;\|\psi_j \|_1/\ell_j^{n+1}
 \quad \text{and}\quad
\|\mu_j\|_{\infty} := \|\psi_j \|_{\infty}/\ell_j,  
\leqno{(2.3)}
$$
where $\ell_j$ denotes the exponent of the test configuration $\mu_j$.
Let $\delta (\mu_j )$ be $\|\mu_j\|_{\infty}/\|\mu_j\|_1$ or $1$
according as $\|\mu_j\|_{\infty} \neq 0$ or $\|\mu_j\|_{\infty} = 0$.
If $\|\mu_j\|_{\infty} \neq 0$, we write $t\in T_{\Bbb R}$ 
as $t = \exp (s/\| \mu_j \|_{\infty})$ for some $s\in \Bbb R$, 
while we require no relation between $s\in \Bbb R$ and $t\in T_{\Bbb R}$ 
if $\|\mu_j\|_{\infty} = 0$.
Note that
$$
\psi_{j}^{\operatorname{SL}}: {T}_{\Bbb R} \to \operatorname{SL}(V_{\ell_j})
$$
is just the restriction of $\psi_j$ to $T_{\Bbb R}$.
Since the group $\operatorname{SL}(V_{\ell_j})$ 
acts naturally on $W_{\ell_j}^*$, 
we can  define a real-valued function $f_{ j}=f_{ j}(s)$ on $\Bbb R$ by
$$
f_{j}(s) \; :=\; \delta (\mu_j )\,\ell_j^{-n}
\log \| \psi_{j}(t)\cdot \hat{X}_{\ell_j} \|^{}_{\operatorname{CH}(\rho^{}_{\ell_j})}, 
\qquad s \in \Bbb R.
\leqno{(2.4)}
$$
Put $\dot{f}_{j}:= df_{j}/ds$.
Here, once $h$ is fixed, the derivative $\dot{f}_{j}(0)$ is bounded from above by a positive constant 
$C$ independent
of the choice of $j$ (see \cite{M}).
Hence we can define $F_{1} (\{\mu_{j}\} )\in \Bbb R \cup \{-\infty \}$ by 
$$
F_{1} (\{\mu_{j}\} )\, :=\; \lim_{s\to -\infty} \{
\varliminf_{j\to \infty}  \dot{f}_{j}(s)\}
\;  \leq\; C,
\leqno{(2.5)}
$$
since the function $\varliminf_{j\to \infty}  \dot{f}_{j}(s)$  is non-decreasing  in $s$ by convexity 
of the function $f_{j}$ (cf. \cite{Zh}; see also \cite{M0}, Theorem 4.5).

\section{Test configurations associated to balanced metrics}

Hereafter, we assume that the polarized algebraic manifold $(X,L)$ is strongly K-stable.
Then by \cite{MN1},
$(X,L)$  is asymptotically Chow-stable, and hence for some
$\ell_0 \gg 1$,
for all $\ell \geq \ell_0$, 
there exists a Hermitian metric $h_{\ell}$ for $L$  such that
$\omega_{\ell} := c_1 (L; h_{\ell})$ is 
a {\it balanced} K\"ahler metric (cf. \cite{D1}, \cite{Zh}) on $(X, L^{\otimes \ell})$
in the sense that
$$
|\sigma^{}_{1}|_{h_{\ell}}^2 +
 |\sigma^{}_{2}|_{h_{\ell}}^2 + \dots +\, |\sigma^{}_{N_{\ell}}|_{h_{\ell}}^2 \; =\; N_{\ell}/c_1(L)^n[X],
\leqno{(3.1)}
$$
where $\{\sigma^{}_{\alpha}\,;\, \alpha = 1,2,\dots,{N_{\ell}}\}$ is an arbitrarily chosen orthonormal basis for 
$(V^{}_{\ell},\rho^{}_{\ell})$. 
Let $\hat{\rho}_{\ell}$ be the associated Hermitian metric on $V_{\ell}$ defined by
$$
\langle\sigma', \sigma'' \rangle^{}_{\hat{\rho}^{}_{\ell}} := \int_X \, (\sigma',\sigma'' )^{}_{h_{\ell}} \, \omega_{\ell}^n,
\qquad \sigma',\sigma'' \in V_{\ell},
$$
where $(\sigma',\sigma'' )_{h_{\ell}}$ denotes the pointwise Hermitian inner product of 
$\sigma$ and $\sigma'$ by the $\ell$-multiple of $h_{\ell}$.
Now  we can find orthonormal bases
$$
\{\sigma^{}_{\ell,1}, \sigma^{}_{\ell,2}, \dots, \sigma^{}_{\ell,N_{\ell}}\}
\;\;\text{ and }\;\; \{\tau^{}_{\ell,1}, \tau^{}_{\ell,2}, \dots, \tau^{}_{\ell,N_{\ell}}\}
$$
for $(V^{}_{\ell},\hat{\rho}^{}_{\ell})$ and $(V^{}_{\ell},\rho^{}_{\ell})$, respectively, such that
$$
\sigma^{}_{\ell,\alpha} \; =\; \lambda^{}_{\ell,\alpha}\tau^{}_{\ell,\alpha},
\qquad \alpha = 1,2,\dots,N_{\ell},
\leqno{(3.2)}
$$
for some positive real numbers $\lambda_{\ell,\alpha}$.
Multiplying $h_{\ell}$ by  
a positive real constant which possibly depends on $\ell$,
we may assume that 
$$
\Pi_{\alpha =1}^{N_{\ell}} \lambda^{}_{\ell,\alpha} = 1.
$$
Then for each $\ell \geq \ell_0$, we have a sequence of points $\hat{\gamma}_k =  (\hat{\gamma}_{k;1}, \hat{\gamma}_{k;2}, \dots ,\hat{\gamma}_{k;N_{\ell}})$, 
$k = 1,2,\dots$, in $\Bbb Q^{N_{\ell}}$ such that $\Sigma_{\alpha = 1}^{N_{\ell}}\hat{\gamma}_{k;\alpha} = 0$
for all $k$, and that
$$
\hat{\gamma}_k \to\; -\,(\log  \lambda^{}_{\ell,1}, \,\log  \lambda^{}_{\ell,2}, \,\dots,\, \log  \lambda^{}_{\ell,N_{\ell}}),
\quad\text{ as $k \to \infty$.}
\leqno{(3.3)}
$$
Let $a_{\ell, k}$ be the smallest positive integer  such that $a_{\ell, k}\hat{\gamma}_k$ is integral. 
By rewriting $a_{\ell,k}\hat{\gamma}_k$ as
${\gamma}_k =  ({\gamma}^{}_{k;1}, {\gamma}^{}_{k;2}, \dots ,{\gamma}^{}_{k;N_{\ell}})$ 
for simplicity,
we now define 
an algebraic group homomorphism $\psi^{}_{\ell, k} : T=\{\,t\in \Bbb C^*\} \to \operatorname{SL}(V_{\ell})$ 
by setting
$$
\psi^{}_{\ell, k} (t)\cdot \tau^{}_{\ell,\alpha} \; :=\; t_{}^{-{\gamma}^{}_{k;\alpha}}\tau^{}_{\ell,\alpha},
\qquad \alpha = 1,2,\dots, N_{\ell},
$$
for all $t \in \Bbb C^*$. 
Let $\{\tau^*_{\ell,\alpha}\,;\, \alpha = 1,2,\dots, N_{\ell}\}$ be the basis for
$V_{\ell}^*$ dual to $\{\tau^{}_{\ell,\alpha}\,;\, \alpha = 1,2,\dots, N_{\ell}\}$ 
defined by
$$
\langle \tau^{}_{\ell,\alpha}, \tau^*_{\ell,\beta}\rangle \; =\;\begin{cases}
&1, \qquad \text{if $\alpha = \beta$},\\
&0, \qquad \text{if $\alpha \neq \beta$.}
\end{cases}
$$
Then $\psi^{}_{\ell, k} (t)\cdot \tau^{*}_{\ell,\alpha} \,=\, 
t_{}^{{\gamma}^{}_{k;\alpha}}\tau^{*}_{\ell,\alpha}$.
Each $\vec{z} = (z_1,z_2, \dots, z_{N_{\ell}}) \in \Bbb C^{N_{\ell}}\setminus \{0\}$
sitting over $(z_1:z_2:\dots :z_{N_{\ell}})\in \Bbb P^{N_{\ell}-1}(\Bbb C )
= \Bbb P^*(V_{\ell}) $ is 
expressible as $\Sigma_{\alpha =1}^{N_{\ell}} z_{\alpha}\tau^{*}_{\ell,\alpha}$, 
and hence the action 
by $t \in \Bbb C^*$ on $\vec{z}$ is written in the form
$$
(z^{}_1, z^{}_2, \dots, z^{}_{N_{\ell}}) \; \mapsto \; (t_{}^{\gamma_{k;1}}z^{}_1, 
t_{}^{\gamma_{k;2}}z^{}_2, \dots , t_{}^{\gamma_{k;N_{\ell}}}z^{}_{N_{\ell}}). 
$$
We now identify $X$ with the subvariety $X_{\ell} := \Phi_{\ell}(X)$ in the projective 
space $\Bbb P^*(V_{\ell}) = \Bbb P^{N_{\ell}-1}(\Bbb C ) = \{(z_1:z_2: \dots :z_{N_{\ell}})\}$ 
via the Kodaira embedding 
$$
\Phi_{\ell} (x ) := (\tau^{}_{\ell,1}(x): \tau^{}_{\ell,2}(x): \dots :\tau^{}_{\ell,2}(x)),
\qquad x \in X.
$$
For each $\ell \geq \ell_0$, we observe that $\operatorname{SL}(V_{\ell})$ acts naturally on $W^*_{\ell}$.
Then by considering the sequence of test configurations 
$$
\mu^{}_{\ell,k} = (\mathcal{X}^{\psi^{}_{\ell, k}},
\mathcal{L}^{\psi^{}_{\ell, k}},{\psi^{}_{\ell, k}}), \qquad k = 1,2,\dots,
$$
associated to $\psi^{}_{\ell, k}$,
we define a real-valued function 
$f^{}_{\ell,k} = f^{}_{\ell,k} (s)$ on the real line $\Bbb R = \{ -\infty < s <+\infty\}$ by
$$
f^{}_{\ell,k} (s) \; :=\; 
\delta (\mu^{}_{\ell,k} ) \,
\ell^{-n}\log  \| \psi_{\ell,k}^{}(t)\cdot \hat{X}_{\ell}\|_{\operatorname{CH}(\rho^{}_{\ell})}.
$$
Here $s\in \Bbb R$ and $t\in \Bbb R_+$ are related by $t = \exp (s/\|\mu_{\ell,k}\|_{\infty})$ 
for $\|\mu_{\ell,k}\|_{\infty}\neq 0$, while we require no relations between $s\in \Bbb R$ and $t\in\Bbb R_+$ 
if $\|\mu_{\ell,k}\|_{\infty} = 0$. 
Put $\dot{f}^{}_{\ell,k}:= df^{}_{\ell,k}/ds$
and $\theta_{s;\ell,k} := (1/2\pi )\log \{(\Sigma_{\alpha = 1}^{N_{\ell}} (n!/\ell^n )\, t^{2{\gamma}^{}_{k;\alpha}} 
|\tau^{}_{\ell,\alpha}|^2)^{1/\ell}\}$. 
Then  on $X_{\ell}$ viewed also as $X$ via $\Phi_{\ell}$, we can write
$$
\psi_{\ell,k}^{}(t)^*(\omega_{\operatorname{FS}}/\ell) 
= \sqrt{-1}\partial\bar{\partial}\theta_{s;\ell,k}, 
\leqno{(3.4)}
$$
where 
$\omega_{\operatorname{FS}} := (\sqrt{-1}/2\pi ) \partial\bar{\partial}\log \{(\Sigma_{\alpha = 1}^{N_{\ell}}(n!/\ell^n )
|z_{\alpha}|^2)_{}^{1/\ell}\}$, and $\psi_{\ell, k}(t)$ is regarded as a mapping from 
$X_{\ell} = (\mathcal{X}^{\psi^{}_{\ell, k}})^{}_1$ to  $\psi_{\ell, k}(t)(X_{\ell}) = (\mathcal{X}^{\psi^{}_{\ell, k}})^{}_t$.
In view of \cite{Zh} (see also \cite{M0} and \cite{S}), we obtain
$$
\dot{f}^{}_{\ell,k} (s) \; =\; \ell\,\delta (\mu_{\ell,k})
\int_X (\partial \theta_{s;\ell,k}/\partial s)\,(\sqrt{-1}\partial\bar{\partial}\theta_{s;\ell,k})^n.
\leqno{(3.5)}
$$
Put
$\nu^{}_{\ell,k} := \|\mu_{\ell,k}\|_{\infty}/a_{\ell,k} 
= \max \{\, |\hat{\gamma}_{k;\alpha}|/\ell\, ;\, \alpha = 1,2,\dots,N_{\ell}\,\}$, where
for the time being, we vary $\ell$ and $k$ independently. Then
$$
(\partial \theta_{s;\ell,k}/\partial s )_{\,|s \,=\, -\nu^{}_{\ell,k}}\; =\;
\frac{\Sigma_{\alpha =1}^{N_{\ell}}\,\hat{\gamma}_{k;\alpha} \,
\exp (-2{\hat{\gamma}^{}_{k;\alpha}})\, |\tau^{}_{\ell,\alpha}|^2}
{\pi \ell \,\nu^{}_{\ell,k} \,\Sigma_{\alpha =1}^{N_{\ell}}\,
\exp (-2{\hat{\gamma}^{}_{k;\alpha}})\,  |\tau^{}_{\ell,\alpha}|^2}.
\leqno{(3.6)}
$$
Now for each integer $r$, let $O(\ell^{\,r})$ denote a function $u$ 
satisfying the inequality $|u| \leq C_0\ell^{\,r}$
for some positive constant $C_0$ independent of the choices of 
$k$, $\ell$, and $\alpha$.  We now fix a positive integer $\ell \gg 1$. Then by (3.3), we obtain
$$
\lambda^{-2}_{\ell,\alpha}\,\exp (-2\hat{\gamma}_{k;\alpha})\, -\,1  \; =\;  O(\ell^{-n-2}), \qquad k\gg 1.
\leqno{(3.7)}
$$
Moreover, in view of  (3.1) and (3.2), the K\"ahler form
$\omega_{\ell}$ is written as
$(\sqrt{-1}/2\pi )\partial\bar{\partial}\log \{ (\Sigma_{\alpha =1}^{N_{\ell}} 
 (n!/\ell^n)\lambda^{\,2}_{\ell,\alpha}|\tau^{}_{\ell,\alpha}|^2)_{}^{1/\ell}\}$. Now by (3.3),
 as $k \to \infty$, we have $\sqrt{-1}\partial\bar{\partial}\theta_{s;\ell,k\,|\,s\,=\,-\nu^{}_{\ell,k}} \to \omega_{\ell}$ in $C^{\infty}$. In particular for $k \gg 1$, we can further assume that
$$
\|
\sqrt{-1}\partial\bar{\partial}\theta_{s;\ell,k\,|\,s\,=\,-\,\nu^{}_{\ell,k}} - \omega_{\ell}\,
\|^{}_{C_{}^m(X)}
\;  = \; O(\ell^{-n-2}),
\leqno{(3.8)}
$$
where we fix an arbitrary integer $m$ satisfying $m \geq 5$.
Hence for each $\ell \gg 1$, we can find a positive integer $k(\ell) \gg 1$ such that
both (3.7) and (3.8) hold for $k  = k(\ell)$. From now on, we assume 
$$
k \;=\; k(\ell ),
\leqno{(3.9)}
$$
and $\nu^{}_{\ell,k }\,=\,\nu^{}_{\ell,k (\ell )}$ will be written as $\nu_{\ell}$ for simplicity. 
Then, since $\ell \nu^{}_{\ell} \geq |\hat{\gamma}^{}_{k;\alpha}|$ for all $\alpha$, 
we have $(\partial \theta_{s;\ell,k}/\partial s )^{}_{|\,s \,=\, -\nu^{}_{\ell}} = O(1)$ by (3.6). Hence 
$$
\int_X (\partial \theta_{s;\ell,k}/\partial s)\,\{ \,(\sqrt{-1}\partial\bar{\partial}\theta_{s;\ell,k})^n \,-\, \omega_{\ell}^{\,n}\,\}^{}_{\,|\,s\,=\,-\nu^{}_{\ell}} \; = \; O(\ell^{-n-2}).
\leqno{(3.10)}
$$
Put $I_1 := \pi \ell \,\nu^{}_{\ell} \,\Sigma_{\alpha =1}^{N_{\ell}}\,
\lambda^2_{\ell,\alpha}\,  |\tau^{}_{\ell,\alpha}|^2 $ and $I_2 := \pi \ell \,\nu^{}_{\ell} \,\Sigma_{\alpha =1}^{N_{\ell}}\,
\exp (-2{\hat{\gamma}^{}_{k;\alpha}})\, |\tau^{}_{\ell,\alpha}|^2 $.
Put also $J_1 := \Sigma_{\alpha =1}^{N_{\ell}}\,\hat{\gamma}_{k;\alpha}
\lambda^2_{\ell,\alpha} |\tau^{}_{\ell,\alpha}|^2 $ and $J_2 :=  \Sigma_{\alpha =1}^{N_{\ell}}
\hat{\gamma}_{k;\alpha}\exp (-2{\hat{\gamma}^{}_{k;\alpha}}) |\tau^{}_{\ell,\alpha}|^2$.
Then by (3.6), we obtain
$$
\int_X \,(\partial \theta_{s;\ell,k}/\partial s )^{}_{\,|s\,=\,-\nu^{}_{\ell}} \;\omega_{\ell}^{\,n}\;
= \; A \,+\, B\, +\,P,
\leqno{(3.11)}
$$
where $ A :=\int_X  \{(J_2/I_2) - (J_2/I_1)\} \,\omega_{\ell}^{\,n}$,
$B := \int_X \{(J_2/I_1) - (J_1/I_1)\}\,\omega_{\ell}^{\,n}$ and 
$P := \int_X \,(J_1/I_1)\,\omega_{\ell}^{\,n}$. Note that $J_2/I_2  = O(1)$
by $\ell\nu^{}_{\ell} \geq |\hat{\gamma}^{}_{k;\alpha}|$,
while by (3.7), $(I_1 - I_2)/I_1 = O(\ell^{-n-2})$. Then 
$$
A \; =\; \int_X\, \frac{J_2}{I_2}\cdot\frac{I_1 - I_2}{I_1}\;\omega_{\ell}^n\; =\; O(\ell^{-n-2}).
\leqno{(3.12)}
$$
On the other hand by (3.7), $J_2 - J_1 = O(\ell^{-n-2}) (\Sigma_{\alpha =1}^{N_{\ell}}\,
|\hat{\gamma}_{k;\alpha}|\,\lambda^2_{\ell,\alpha}\,    |\tau^{}_{\ell,\alpha}|^2)$.
From this together with $\ell\nu^{}_{\ell} \geq |\hat{\gamma}^{}_{k;\alpha}|$,
we obtain 
$$
B \; =\; \int_X \frac{J_2 - J_1}{I_1} \;\omega_{\ell}^{\,n} \; = \; O(\ell^{-n-2}).
\leqno{(3.13)}
$$
By (3.2),  $ I_1 = \pi \ell \,\nu^{}_{\ell} \,\Sigma_{\alpha =1}^{N_{\ell}}\,
  |\sigma^{}_{\ell,\alpha}|^2$ and $J_1 := \,\Sigma_{\alpha =1}^{N_{\ell}}\,\hat{\gamma}_{k;\alpha}
\,  |\sigma^{}_{\ell,\alpha}|^2 $.  
Note also that $a^{}_0 :=  \delta (\mu^{}_{\ell,k} )$ satisfies $ 0 < a^{}_0 \leq \ell^n$.
Put $a_1 := c_1(L)^n[X]$.
In view of (3.1) and (3.5), by adding up (3.10), (3.11), (3.12) and (3.13), we obtain
$$
\begin{cases}
&\dot{f}^{}_{\ell,k} (-\nu^{}_{\ell}) \, =\, 
\ell\,a_0\,\int_X \{\,(\partial \theta_{s;\ell,k}/\partial s)\,(\sqrt{-1}\partial\bar{\partial}\theta_{s;\ell,k})^n\}^{}_{\,|s\,=\,-\nu^{}_{\ell}}\\
&=\;  a^{}_0\,\{\, \ell\,P  \,+\,O(\ell^{-n-1})\,\}
=\,\int_X \frac{a^{}_0\Sigma_{\alpha =1}^{N_{\ell}}\,\hat{\gamma}_{k;\alpha}
  |\sigma^{}_{\ell,\alpha}|_{h_{\ell}}^2\;}{\;\pi \,\nu^{}_{\ell} \,\Sigma_{\alpha =1}^{N_{\ell}}
  |\sigma^{}_{\ell,\alpha}|_{h_{\ell}}^2\;}\, \omega_{\ell}^n \,+\, O(\ell^{-1})\\
  & =\; 
 a_0\, a_1\,(\Sigma_{\alpha =1}^{N_{\ell}}\,\hat{\gamma}_{k;\alpha})\,
 (\pi\, \nu^{}_{\ell}\, N_{\ell})^{-1} \, +\, O(\ell^{-1}) \; =\;O(\ell^{-1}) , 
\end{cases}
\leqno{(3.14)}
$$
where in the last line, we used the equality $\Sigma_{\alpha =1}^{N_{\ell}}\,\hat{\gamma}_{k;\alpha} = 0$.
In the next section, the sequence of test configurations $\mu^{}_{\ell,k(\ell )} = (\mathcal{X}^{\psi^{}_{\ell, k(\ell )}},
\mathcal{L}^{\psi^{}_{\ell, k(\ell )}},{\psi^{}_{\ell, k(\ell )}})$, $\ell \geq \ell_0$, for $(X,L)$ will be considered.

\section{Proof of Main Thorem}

In this section, under the same assumption as in the previous section, we shall show that 
$c_1(L)$ admits a constant scalar curvature K\"ahler metric. 
Put 
$$
\nu^{}_{\infty}\;  :=\; \sup_{\ell} \; \nu^{}_{\ell},
$$
where the supremum is taken over all positive integers $\ell$ satisfying $\ell \geq \ell_0$.
Then the following cases are possible:

\medskip
 \qquad Case 1:  \; $\nu^{}_{\infty}\, =\, +\infty$.
 \qquad Case 2:  \; $\nu^{}_{\infty} \,<\, +\infty$.
 
 \medskip\noindent
{\em Step \/$1$.}  If Case 1 occurs, then an increasing subsequence $\{\,\ell_j\,; \,j=1,2,\dots \,\}$ 
of  $\{\,\ell\in \Bbb Z\,;\, \ell \geq \ell_0\,\}$ 
can be chosen in such a way that $\{\nu^{}_{\ell_j}\}$ is a  monotone increasing sequence 
satisfying
$$
\lim_{j\to \infty} \,\nu^{}_{\ell_j}\; =\; +\infty.
\leqno{(4.1)}
$$
For simplicity, the functions $f^{}_{\ell_j, k(\ell_j )}$ will be written as $f_{j}$, while  
we write the test configurations 
$$
\mu^{}_{\ell_j, k(\ell_j )}= (\mathcal{X}^{\psi^{}_{\ell_j, k(\ell_j )}},
\mathcal{L}^{\psi^{}_{\ell_j, k(\ell_j )}},{\psi^{}_{\ell_j, k(\ell_j )}}), \qquad j=1,2,\dots,
$$
as $\mu_j = (\mathcal{X}_j, \mathcal{L}_j, \psi_j )$.
Now by (3.14), there exists a positive constant $C$ independent of $j$
 such that 
 $$
-\, C /\ell_j\;\leq \;  \dot{f}_{j}^{}(-\nu^{}_{\ell_j})
 $$
for all $j$. On the other hand, for all positive integers $ j' $ satisfying $j' \geq j$, we have $-\nu^{}_{\ell_{j'}}
\leq -\nu^{}_{\ell_j}$ by monotonicity. Since the function $ \dot{f}_{j'}^{}(s)$ in $s$ is non-decreasing, 
we obtain
$$
 -\, C/\ell_{j'} \;\leq \; \dot{f}_{j'}^{}(-\nu^{}_{\ell_{j'}})\; \leq \;  \dot{f}_{j'}^{}(-\nu^{}_{\ell_{j}}).
 \leqno{(4.2)}
$$
We here observe that $ -\, C/\ell_{j'}  \to 0$ as $j' \to \infty$.  It now follows from (4.2) that, for each fixed $j$, 
$$
\varliminf_{j' \to \infty}\dot{f}_{j'}^{}(-\nu^{}_{\ell_{j}}) \; \geq\; 0.
$$
Since the function  $\varliminf_{j' \to \infty}\dot{f}_{j'}^{}(s)$ in $s$ is non-decreasing, 
we therefore obtain 
$$
\varliminf_{j' \to \infty}\dot{f}_{j'}^{}(s ) \; \geq \; 0
\;\, \,\text{for all \,$ s \geq -\nu^{}_{\ell_{j}}$,}
$$
while this holds for all positive integers $j$. Then by (4.1), 
$\varliminf_{j' \to \infty}\dot{f}_{j'}^{}(s )$ is a nonnegative function in $s$ on the whole 
real line $\Bbb R$.
Hence 
$$
F_1( \{\mu_j\}) \; =\; \lim_{s\to -\infty}\{\varliminf_{j' \to \infty}\dot{f}_{j'}^{}(s ) \} \; \geq \; 0.
$$
Now by the strong K-stability of $(X,L)$, we obtain $F_1(\{\mu_j \}) = 0$, so that
$\mu_j$ are trivial  for all $j\gg 1$. Then $\psi_{\ell_j, k(\ell_j)}$ are trivial 
for all $j\gg 1$. This usually gives us a contradiction. Even if not, however, by assuming the triviality 
of $\mu_j$ for all $j\gg 1$, we proceed as follow.
By (3.4), for all $s \in \Bbb R$, we obtain
$$
\sqrt{-1} \partial\bar{\partial}\theta_{s;\ell_j, k(\ell_j )}\; =\;
({\omega_{\operatorname{FS}}}/\ell_j )_{|X_{\ell_j}} \; =\; \Phi_{\ell_j}^*({\omega_{\operatorname{FS}}}/\ell_j ),
\; \qquad   j \gg 1,
$$
 by identfying $X_{\ell_j}$ with $X$ via $\Phi_{\ell_j}$, 
where by \cite{Ze}, 
 $ \|\Phi_{\ell_j}^*({\omega_{\operatorname{FS}}}/\ell_j )\, -\, \omega\,\|^{}_{C^5(X)}\,=\,O(\ell_j^{-2})$. 
From this together with  (3.8), we obtain
$$
\|\,\omega \,-\,\omega^{}_{\ell_j}\,\|^{}_{C_{}^m(X)} \; =\; O(\ell_j^{-2}), \qquad j \gg 1.
\leqno{(4.3)}
$$
Let $S_{\omega}$ be the scalar curvature function for $\omega$. Then by \cite{Lu} (see also \cite{Ze}),
we obtain the following asymptotic expansion:
$$
1+ (S_{\omega}/2)\ell_j^{-1} + O(\ell_j^{-2})\; =\; 
\Sigma_{\alpha =1}^{N_{\ell_j}}
(n!/\ell_j^n ) \,|\tau^{}_{\ell_j, \alpha}|^2_{h} \; =\; B_{\ell_j}(\omega ),
\leqno{(4.4)}
$$
where for every K\"ahler form $\theta$ in $c_1(L)_{\Bbb R}$,
$B_{\ell_j}(\theta )$ denotes the $\ell_j$-th asymptotic Bergman kernel for $(X,\theta )$.
On the other hand, for $\ell \gg 1$,  we observe that $N_{\ell}$ is a polynomial in $\ell$. Since each $\omega_{\ell_j}$ is balanced, by setting $\ell = \ell_j$ in (3.1) and 
dividing both sides of the equality by $\ell_j^n /n!$, we obtain (cf. \cite{MA}, (1.4)) 
$$
1\, + \,C^{}_0\, \ell^{-1}_j +\,O( \ell_j^{-2})\; =\; \Sigma_{\alpha =1}^{N_{\ell_j}}
(n!/\ell_j^n ) \,|\sigma^{}_{\ell_j, \alpha}|^2_{h^{}_{\ell_j}}
\; =\; B_{\ell_j}(\omega_{\ell_j} ),
\leqno{(4.5)}
$$
where $C^{}_0$ is a real constant independent of the choice of $j$. 
In view of (4.3), by comparing (4.4) with (4.5), we now conclude that $S_{\omega}/2 = C_0$. Hence $\omega$ is a constant scalar curvature K\"ahler metric in the class $c_1(L)_{\Bbb R}$.

\medskip\noindent
{\em Step \/$2$.} Suppose that Case 2 occurs. Put $\hat{\lambda}_{\ell, \alpha} := -(1/\ell )\log \lambda_{\ell, \alpha}$. 
Then by (3.3), we may assume that $k = k(\ell )$ in (3.9) is chosen in such a way that 
$$
\hat{\gamma}_{k (\ell );\alpha} -1\, \leq \; \ell\, \hat{\lambda}_{\ell, \alpha}\; \leq \,\hat{\gamma}_{k(\ell );\alpha} +1,
\qquad \alpha = 1,2,\dots, N_{\ell},
\leqno{(4.6)}
$$
for all $\ell$ with $\ell \geq \ell_0$. Then for each $\ell$, by using the notation in Definition 5.3 in Appendix, 
we have an $\ell$-th root 
$$
(\mathcal{Y}^{(\ell )}, \mathcal{Q}^{(\ell )}, D^{(\ell )}, \varphi^{}_{\ell} ),
\qquad \ell \geq \ell_0,
$$ 
of the test configuration $\mu_{\ell,k(\ell )}$ in Section 3. Let $\chi^{}_{\ell,\beta}$, $\beta =1,2,\dots, N_1$, be the weights of the $T_{\Bbb R}$-action via $\varphi_{\ell}^{\operatorname{SL}}$ on $V_1^*$, where
$V_1:= H^0(X,L)$.
Put $\hat{\chi}^{}_{\ell,\beta} := \chi^{}_{\ell, \beta}/a_{\ell, k(\ell )}$.  
For $\ell$ with $\ell \geq \ell_0$,
let $\alpha$ and $\beta$ be arbitrary integers satisfying
$1\leq \beta \leq N_1$ and $1\leq \alpha \leq N_{\ell}$.
By (4.6) together with the definition of $\nu_{\ell,k}$,
we easily see from the inequality $\nu_{\infty}  < +\infty$ that
$$
|\hat{\lambda}_{\ell, \alpha}| \leq C_1\quad \text{ and }\quad
|\hat{\chi}^{}_{\ell,\beta} |\leq C_1,
\leqno{(4.7)}
$$
where $C_1$ is a positive real constant independent of 
the choices of $\ell$, $\alpha$ and  $\beta$ (see \cite{MN2} for the second inequality of (4.7)\,; see also \cite{M2}).
Let $Z_{\ell} := (\varphi_{\ell})_* (t\partial/\partial t) \in \frak{sl}(V_1)$ be the infinitesimal generator
for the one-parameter group $\varphi_{\ell}^{\operatorname{SL}}$. 
Then by setting $\hat{Z}_{\ell}:= Z_{\ell}/a_{\ell,k(\ell)}$, 
we obtain
$$
\hat{Z}_{\ell} \cdot \kappa^{}_{\ell,\beta }\; =\; -\,\hat{\chi}^{}_{\ell,\beta}\kappa^{}_{\ell,\beta},
\qquad  \beta = 1,2,\dots,N_1,
$$
for a suitable orthonormal basis $\{ \kappa_{\ell, 1}, \kappa_{\ell, 2}, 
\dots, \kappa_{\ell, N_1}\}$  for $(V_1,\rho_1) $.
For the sequence $\{\,\hat{Z}_{\ell}\,;\,\ell \geq \ell_0\,\}$, by choosing its suitable
subsequence  
$$
\{\,\hat{Z}_{\ell_j}\,;\,j=1,2,\dots\,\},
$$
we obtain
real numbers $\hat{\chi}^{}_{\infty,\beta}\in \Bbb R$, $\beta = 1,2,\dots,N_1$, and an orthonormal
basis $\{ \kappa^{}_{\infty, 1}, \kappa^{}_{\infty, 2}, 
\dots, \kappa^{}_{\infty, N_1}\}$ for $V_1 $ such that, for all $\beta$,
$$
\kappa^{}_{\ell_j,\beta} \to \kappa^{}_{\infty,\beta} \quad\text{ and }\quad
\hat{\chi}^{}_{\ell_j,\beta} \to \hat{\chi}^{}_{\infty,\beta},
$$ 
as $j \to \infty$. Hence we can define  $\hat{Z}_{\infty} \in \frak{sl}(V_1)$ such that
$\hat{Z}_{\infty} \cdot \kappa^{}_{\infty,\beta }\, =\, -\,\hat{\chi}^{}_{\infty,\beta}\kappa^{}_{\infty,\beta}$
for all $\beta$. Then we have the following convergence in $C^{\infty}$:
$$
\hat{Z}_{\ell_j} \to \hat{Z}_{\infty}, 
\qquad \text{ as $j \to \infty$.}
\leqno{(4.8)}
$$
For each $\ell$, in view of the relation $t = \exp (s/\|\mu_{\ell, k(\ell )}\|_{\infty}) $, 
$s = -\nu_{\ell}$ corresponds to $t = \hat{t}_{\ell}$, where 
$\hat{t}_{\ell} := \exp (-\nu_{\ell}/\|\mu_{\ell, k(\ell )}\|_{\infty}) = \exp (-1/a_{\ell, k(\ell )})$. 
Until the end of this section, test configurations $\mu^{}_{\ell,k(\ell )} $ for $(X,L)$ will be written simply as
$$
\mu_{\ell} = (\mathcal{X}^{(\ell )},
\mathcal{L}^{(\ell )},\psi^{}_{\ell}),  \qquad \ell \geq \ell_0.
$$
For the test configuration $\mu_{\ell}$,
each $t \in T$ not as a complex number but as an element of the group $T$ of transformations on $\mu_{\ell}$
will be written as $g_{\mu_{\ell}}(t)$.
For the Kodaira embedding $\Phi_{\ell}: X  \hookrightarrow \Bbb P^{N_{\ell}-1}(\Bbb C ) $ 
in Section 3, we consider $\Bbb C^{N_{\ell}} \setminus \{0\} = \{(z_1,z_2, \dots , z_{N_{\ell}})\neq 0\}$ 
over $\Bbb P^{N_{\ell}}(\Bbb C )$,
so that $z= (z_1,z_2,\dots, z_{N_{\ell}})$ sits over 
$[z] = (z_1:z_2: \dots :z_{N_{\ell}})$. 
 Since the restriction of $\mathcal{O}_{\Bbb P^{N_{\ell}-1}}(\Bbb C )$
to ${X_{\ell}}$ 
is viewed as $L$ by identifying $X$ with its image $X_{\ell} := \Phi_{\ell}(X)$, we can write
$$
z_{\alpha} {}_{|X_{\ell}} \; =\; \tau_{\ell,\alpha}, \qquad \alpha = 1,2,\dots, N_{\ell},
$$
for the orthonormal basis
$\{\tau_{\ell,1},\tau_{\ell,2}, \dots, \tau_{\ell,N_{\ell}}\}$  of 
$(V_{\ell}, \rho_{\ell})$.
We now define a Hermitian metric  $\phi_{\ell}$ for $L^{-1}$
by setting, for all $[z] = \Phi_{\ell}(x)$ in $X_{\ell}$,
$$
\phi^{}_{\ell} ([z])\,:=\, \{\,(n!/\ell^n)\Sigma_{\alpha =1}^{N_{\ell}}|z^{}_{\alpha}|^2\,\}^{1/\ell}
\, =\, \{\,(n!/\ell^n)\Sigma_{\alpha =1}^{N_{\ell}}|\tau^{}_{\ell,\alpha}(x)|^2\,\}^{1/\ell},
$$ 
where the line bundle $L^{-\ell}$ on $X$ is viewed as
 the dual 
 $\{\mathcal{L}^{(\ell )}{}_{|X_{\ell}}\}_{}^{-1}$ 
 of the line bundle ${\mathcal{L}^{(\ell )}}$ 
restricted to $\mathcal{X}_1^{(\ell )} \, (=X_{\ell})$. Let $\mathcal{K}_t$, $t \neq 0$,  denote the set of all Hermitian metrics on the 
line bundle $\{\mathcal{L}^{(\ell )}{}_{|\mathcal{X}_t^{(\ell )}}\}_{}^{-1}$.
Then the action by $g_{\mu_{\ell}}(t) $ takes 
$\mathcal{K}_1$ 
to $\mathcal{K}_t$. For instance, $g_{\mu_{\ell}}(t)$ takes the point $z= (z_1, z_2, \dots , z_{N_{\ell}})$
to $g_{\mu_{\ell}}(t)\cdot z = 
(t^{\gamma_{k(\ell ),1}}z_1, t^{\gamma_{k(\ell ),2}}z_2, \dots, t^{\gamma_{k(\ell ),N_{\ell}}}z_{N_{\ell}})
$, while for each $[z]\in X_{\ell}$, $\phi_{\ell}([z])$ is mapped to
the point $g_{\mu_{\ell}}(t)\cdot \phi_{\ell}([z])$ defined by
$$
\{\,(n!/\ell^n)\Sigma_{\alpha =1}^{N_{\ell}}|g_{\mu_{\ell}}(t)\cdot z_{\alpha}|^2\,\}^{1/\ell}
\;=\; \{\,(n!/\ell^n)\Sigma_{\alpha =1}^{N_{\ell}}|t|^{2\gamma_{k(\ell ),\alpha}}| z_{\alpha}|^2\,\}^{1/\ell},
$$
and this defines $g_{\mu_{\ell}}(t)\cdot \phi_{\ell}\in \mathcal{K}_t$.
Now by \cite{Ze},
$u_{\ell}:= (1/2\pi) \log (\phi_{\ell}/h^*)$ viewed as a function on $X$ can be estimated in the form
$$
\|u_{\ell}\|_{C^{m+2}_{}(X)} = O(\ell ^{-2}),
\leqno{(4.9)}
$$
where the dual $h^*$ of $h$ is viewed as a Hermitian metric for the line bundle 
$L^{-1}$.
Put $\omega (\ell, t) :=  (\sqrt{-1}/2\pi )\partial\bar{\partial}
\log (g_{\mu_{\ell}}(t)^*\{g_{\mu_{\ell}}(t)\cdot h^*\})$,
$t \neq 0$.
For the Fubini-Study form $\omega_{\operatorname{FS}}$ in Section 3, 
its restriction to $X_{\ell} \, (=X)$ is written as
$$
\omega_{\operatorname{FS}}\,{}_{|X_{\ell}}\;=\; (\sqrt{-1}\,\ell/2\pi )\partial\bar{\partial}\log \phi_{\ell}.
$$
Since $\psi_{\ell}(t)^*(\omega_{\operatorname{FS}}/\ell )$
$= (\sqrt{-1}/2\pi )\partial\bar{\partial}\log (g_{\mu_{\ell}}(t)^*\{g_{\mu_{\ell}}(t)\cdot \phi_{\ell}\})$
on $X_{\ell}$ (see (3.4)),
we can rewrite it in the form (see \cite{MR})
$$
\psi_{\ell}(t)^*(\omega_{\operatorname{FS}}/\ell ) \,{}_{|X_{\ell}}\;=\; 
 \omega (\ell, t )\,+\,
\sqrt{-1}\,\partial\bar{\partial}u_{\ell}.
\leqno{(4.10)}
$$
Let us consider the test configuration  $\bar{\mu}_{\ell} := (\mathcal{Y}^{(\ell )}, \mathcal{Q}^{(\ell )}, \varphi_{\ell})$ 
for $(X,L)$ of exponent $1$. 
Each $t \in T$, not as a complex number but as an element of the group $T$ of transformations
on $\bar{\mu}_{\ell}$, will be denoted by $g_{\bar{\mu}_{\ell}}(t)$. Then by (5.8) in Appendix, 
we also have the expression
$$
\omega (\ell, t) \;=  \; (\sqrt{-1}/2\pi )\partial\bar{\partial}
\log (g_{\bar{\mu}_{\ell}}(t)^*\{g_{\bar{\mu}_{\ell}}(t)\cdot h^*\}),
\;\;\quad t\in T_{\Bbb R},
\leqno{(4.11)}
$$
since for each such $t$, 
the action of $g_{\mu_{\ell}}(t)$ on $|\mathcal{L}|^{2/\ell}$ coincides with the action of 
$g_{\bar{\mu}_{\ell}}(t)$ 
on $|\mathcal{Q}|^2$ up to constant scalar multiplication, where constant scalar multiplication arises from the action on the factor $|\zeta |^{2/\ell}$. Since $\omega (\ell, t)$ doesn't change even if 
$g_{\bar{\mu}_{\ell}}(t)\cdot h^*$ in (4.11)
is replaced by $C(t)g_{\bar{\mu}_{\ell}}(t)\cdot h^*$
for a positive real constant $C(t)$ possibly depending on $t$.  
Hence we may consider the action by  $g_{\bar{\mu}_{\ell}}(t)$ on $|\mathcal{L}|^{2/\ell}$ modulo constant scalar multiplication.
In this sense, for each $t \in T_{\Bbb R}$, the action by $g_{\bar{\mu}_{\ell}}(t)$ in (4.11)  is induced by 
the action by the element $\varphi_{\ell}^{\operatorname{SL}}(t)$ in $\operatorname{SL}(V_1)$.
In particular for $t =  \hat{t}_{\ell}$, 
$$
\text{$g_{\bar{\mu}_{\ell}}(\hat{t}_{\ell})\,$'s action is induced by  $\,\varphi_{\ell}^{\operatorname{SL}}(\hat{t}_{\ell})
\; =\;  \exp (-\hat{Z}_{\ell})\,\in \,\operatorname{SL}(V_1)$.}
\leqno{(4.12)}
$$
For $\theta_{s;\ell,k(\ell )} := (1/2\pi )\log \{(\Sigma_{\alpha = 1}^{N_{\ell}} (n!/\ell^n )\, t^{2{\gamma}^{}_{k(\ell );\alpha}} 
|\tau^{}_{\ell,\alpha}|^2)^{1/\ell}\}$
in Section 3,  at the point $s = -\nu_{\ell}$, we see from (3.4) that
$$
\sqrt{-1}\partial\bar{\partial}\theta_{s;\ell,k(\ell )}\,{}_{|s = -\nu_{\ell} }
\; =\;\psi_{\ell}(\hat{t}_{\ell})^*(\omega_{\operatorname{FS}}/\ell ) \,{}_{|X_{\ell}}.
\leqno{(4.13)}
$$
Then by (3.8), (4.10) and (4.13), 
$$
\| \omega_{\ell} - \omega (\ell, \hat{t}_{\ell}) -\sqrt{-1}\partial\bar{\partial}u_{\ell}\|_{C^m_{}(X)} \; =\; O(\ell^{-n-2}).
\leqno{(4.14)}
$$
For the element $\hat{Z}_{\infty}$ of $\frak{sl}(V_1)$ in (4.8), 
we now define a subset $\mathcal{Y}_{\Bbb R}^{(\infty )}$ of 
$\Bbb R \times \Bbb P^*(V_1)$ as the closure of 
$$
\bigcup_{s\in \Bbb R}\; \{\pm \exp s\} \times \exp (s\hat{Z}_{\infty})(X_1)
$$
in the real manifold $\Bbb R \times \Bbb P^*(V_1)$, where $X_1$ is the image $\Phi_1 (X)$  of $X$ under the Kodaira embedding 
$$
\Phi_1: X \to \Bbb P^*(V_1)
$$ 
associated to the complete linear system $|L|$ on $X$.
By the projection of $\Bbb R \times \Bbb P^*(V_1)$ to the first factor $\Bbb R$, 
we see that $\mathcal{Y}_{\Bbb R}^{(\infty )}$ has a natural structure of a fiber space over $\Bbb R$.
Let $\mathcal{Q}^{(\infty )}$ denote the restriction to 
$\mathcal{Y}_{\Bbb R}^{(\infty )}$ of the pullback $\operatorname{pr}_2^*\mathcal{O}_{\Bbb P^*(V_1)}(1)$,
where $\operatorname{pr}_2: \Bbb R \times \Bbb P^*(V_1) \to \Bbb P^*(V_1)$ is the projection 
to the second factor. Then the $T_{\Bbb  R}$-action on $\mathcal{Y}_{\Bbb R}^{(\infty )} $ induced by
$$
T_{\Bbb R} \times (\Bbb R \times \Bbb P^*(V_1)) \to \Bbb R \times \Bbb P^*(V_1),
\;\; 
(\exp s, (r,x)) \mapsto ( \,(\exp s )r, \exp (s\hat{Z}_{\infty})\cdot x\,),
$$
naturally lifts to 
a  $T_{\Bbb  R}$-action on $\mathcal{Q}^{(\infty )}$.
This action is induced by the  Lie group homomorphism 
$\varphi_{\infty}: \Bbb R_+ \to \operatorname{SL}(V_1)$ 
defined by 
$$
\varphi_{\infty}(t) := \exp (\,(\log t) \hat{Z}_{\infty}\,),
\qquad t \in \Bbb R_+.
$$
For $\bar{\mu}_{\infty} := (\mathcal{Y}_{\Bbb R}^{(\infty )}, \mathcal{Q}^{(\infty )},\varphi_{\infty})$, 
each $t \in T_{\Bbb R}$ not as a real number but as an element of the  group $T_{\Bbb R}$ of transformations
on $\bar{\mu}_{\infty}$ will be written as $g_{\bar{\mu}_{\infty}}(t)$. 
Consider the action by  $g_{\bar{\mu}_{\infty}}(\hat{t}^{}_{\ell})$ on $| \mathcal{Q}^{(\infty )}|^2$ 
modulo constant scalar multiplication.
For $\hat{t}_{\infty} := 1/e$,
we have $\varphi_{\infty}(\hat{t}^{}_{\infty})\; =\; \exp (- \hat{Z}_{\infty})$, and hence
$$
\text{$g_{\bar{\mu}_{\infty}}(\hat{t}^{}_{\infty})\,$'s action is induced by  $\,\varphi_{\infty}(\hat{t}^{}_{\infty})
\; =\;  \exp (-\hat{Z}_{\infty})\,\in \,\operatorname{SL}(V_1)$.}
\leqno{(4.15)}
$$
Put $\omega_{\infty} := (\sqrt{-1}/2\pi )\partial\bar{\partial}
\log (g_{\bar{\mu}_{\infty}}(\hat{t}^{}_{\infty})^*\{g_{\bar{\mu}_{\infty}}(\hat{t}^{}_{\infty})\cdot h^*\})$ (cf.~Remark 5.9).
By (4.11), (4.12) and (4.15), 
it follows from (4.8) that
$$
\omega (\ell_j, \hat{t}_{\ell_j}) \to \omega_{\infty} \text{\, in\, $C^{\infty}$},
\qquad \text{ as $j \to \infty$.}
\leqno{(4.16)}
$$
Then by (4.9), (4.14) and (4.16), 
$$
\omega_{\ell_j} \to \omega_{\infty} \,\text{ in $C^m_{}$}, \quad \text{as $j \to \infty$.}
\leqno{(4.17)}
$$
By (4.17), given a sufficiently small $\varepsilon >0$, there exists a $j_0 \gg 1$
such that
$\| S_{\omega_{\ell_j}} - S_{\omega_{\infty}} \|_{C^0(X)} \leq \varepsilon$
for all $ j \geq j_0$.
Hence by  \cite{Lu} (see also  \cite{Ze}),
$$
|\,\ell_j \{B_{\ell_j}(\omega_{\ell_j}) - \hat{N}_{\ell_j}\} - \ell_j \{B_{\ell_j}(\omega_{\infty}) - \hat{N}_{\ell_j}\}\,|\;
\leq\; \varepsilon /2 \,+\, O(1/\ell_j ), 
\quad j\geq j_0,
$$
where $\hat{N}_{\ell_j}:= (n!/\ell_j^n)N_{\ell_j}/c_1(L)^n[X]$.
On the other hand, since each $\omega_{\ell_j}$ is balanced, 
we have $B_{\ell_j}(\omega_{\ell_j}) = \hat{N}_{\ell_j}$  for all $j$. 
It then follows that
$$
|\,\ell_j \{B_{\ell_j}(\omega_{\infty}) - \hat{N}_{\ell_j}\}\,| \; \leq\; \varepsilon /2 \,+\, O(1/\ell_j ), 
\quad j\geq j_0.
$$
Hence, since $\hat{N}_{\ell_j} = 1 +  C_0 \ell_j^{-1}+ O(\ell_j^{-2})$ for a real constant $C_0$ independent of 
$j$,
again by \cite{Lu} (see also \cite{Ze}) applied to $\omega_{\infty}$, we obtain 
$$
|\, (S_{\omega_{\infty}}/2)\, -\, C_0\, |\; \leq \; \varepsilon/2  \, + \, O(1/\ell_j ), 
\qquad j\geq j_0,
$$
so that by letting $j \to \infty$, we have $|\, (S_{\omega_{\infty}}/2)\, -\, C_0\, |\, \leq \, \varepsilon/2$.
Since $\varepsilon >0$ can be chosen arbitrarily, we obtain $S_{\omega_{\infty}} = 2 C_0$, as required.

\medskip\noindent
{\em Remark \/$4.18$.}
The $(1,1)$-form $\omega_{\infty}$ on $X$ is positive-definite as follows:
For each $t \in T_{\Bbb R}$  viewed as a real number, the fiber of $\mathcal{Y}_{\Bbb R}^{(\infty )}$ over 
$t \in \Bbb R\setminus \{0\}$ will be denoted by $\mathcal{Y}_t$,
where $\mathcal{Y}_t \cong X$ biholomorphically.
For simplicity, the fiber $(\mathcal{Q}^{(\infty )})^{}_t$ of $\mathcal{Q}^{(\infty )}$ over $t$ will be written as 
$\mathcal{Q}_t$,
and $g_{\bar{\mu}_{\infty}}(\hat{t}_{\infty}^{})$ will be written as $g$.
Then  $ g$ takes 
$\mathcal{Y}_1$ 
holomorphically onto $\mathcal{Y}_{\hat{t}^{}_{\infty}}$. Hence 
$$
\omega_{\infty} = (\sqrt{-1}/2\pi )g^*\partial\bar{\partial}
\log (g\cdot h^*).
\leqno{(4.19)}
$$
Moreover $ g:
\mathcal{Y}_1 \to \mathcal{Y}_{t'}$ lifts
holomorphically to a map, denoted also by $ g$ by abuse of terminology,
 of $\mathcal{Q}_1$ 
onto $\mathcal{Q}_{\hat{t}_{\infty}^{}}$. 
By choosing a local base $b$ for $\mathcal{Q}_1$ 
on an open subset $U$ of $\mathcal{Y}_1$, we can write $h^*$  as
$Hb\bar{b}$ 
for some positive real-valued function $H$ on $U$. 
Since $\omega = c_1(L; h)$ is K\"ahler,
$\sqrt{-1}\partial\bar{\partial}\log H$ is positive-definite on $U$.
Then by $g\cdot h^* = \,
(H\circ g^{-1})\,g(b) \overline{g(b)}$, we see that
$$
\sqrt{-1}\partial\bar{\partial}\log (g\cdot h^*)\; =\;\sqrt{-1}\partial\bar{\partial}\log (H\circ g^{-1})
$$
is positive-definite on $g(U)$. From this together with (4.19),  we now conclude that  $\omega_{\infty} $ is 
positive-definite.


\section{Appendix}

In this appendix, we consider a test configuration $\mu = (\mathcal{X},\mathcal{L},\psi )$ for $(X,L)$, 
and let $\pi: \mathcal{X}\to \Bbb A^1$ be the associated $T$-equivariant projective morphism. 
For the exponent $\ell$ of $\mu$, $\psi$ is an algebraic group homomorphism
$$
\psi : \Bbb C^* \to \operatorname{GL}(V_{\ell}),
$$
and by choosing a Hermitian metric $h$ for $L$,
we endow $V_{\ell}:= H^0(X,L^{\otimes \ell})$ with the Hermitian metric 
$\rho_{\ell}$  as in the introduction.

\medskip\noindent
{\em Definition \/$5.1$.}  A pair $(\hat{\mathcal{X}},\hat{\mathcal{L}})$ of a non-singular irreducible algebraic variety $\hat{\mathcal{X}}$ and an invertible sheaf
 $\hat{\mathcal{L}}$ over $\hat{\mathcal{X}}$  is called a 
 {\it $T$-equivariant desingularization} \/of $(\mathcal{X},\mathcal{L})$, 
 if there exists a $T$-equivariant proper birational morphism $\iota  :  \hat{\mathcal{X}} \to \mathcal{X}$, 
 isomorphic over $\mathcal{X}\setminus \mathcal{X}_0$, such that $\hat{\mathcal{L}} = \iota^*\mathcal{ L}$.
 

\medskip\noindent
{\bf Theorem 5.2.} 
{\em There exist a $T$-equivariant desingularization $(\hat{\mathcal{X}},\hat{\mathcal{L}})$ of $(\mathcal{X},\mathcal{L})$
and a test configuration $(\mathcal{Y},\mathcal{Q},\varphi )$ for $(X,L)$ of exponent $1$ 
such that
$$
\hat{\mathcal{L}} \, =\, \mathcal{O}_{\hat{\mathcal{X}}}(\hat{D})\otimes \eta^*\mathcal{Q}^{\otimes\ell},
$$
where $\eta : \hat{\mathcal{X}} \to \mathcal{Y}$ is a $T$-equivariant proper birational morphism,
isomorphic over $\mathcal{Y}\setminus \mathcal{Y}_0$,
 and $\hat{D}$ is a divisor on $\hat{\mathcal{X}}$ sitting in 
$\hat{\mathcal{X}}_0$ set-theoretically. 
}

\medskip\noindent
{\em Definition \/$5.3$.} Taking the $\Bbb Q$-divisor $D:= \hat{D}/\ell$ on $\hat{\mathcal{X}}$, we call the quadruple $(\mathcal{Y},\mathcal{Q}, D, \varphi )$ an {\it $\ell$-th root\/} of the test configuration $(\mathcal{X},\mathcal{L},\psi )$.


\medskip\noindent
{\em Proof}:
Consider the relative Kodaira embedding
$$
\mathcal{X} \;\hookrightarrow \; \Bbb A^1 \times \Bbb P^* (V_{\ell})
$$
whose restriction $\mathcal{X}_z \hookrightarrow \{z\}\times \Bbb P^*(V_{\ell})$ over each $z \in \Bbb A^1$ is the Kodaira embedding of 
$\mathcal{X}_z$ by the complete linear system $|\mathcal{L}_z|$. 
Let $H$ be a general member in the complete linear system $|L|$ for the line bundle $L$ on $X$.
By the identification $X = \mathcal{X}_1$, we view $H$ as a divisor on $\mathcal{X}_1$.
Then on the projective bundle $\Bbb A^1 \times \Bbb P^* (V_{\ell})$, a $T$-invariant irreducible reduced divisor 
$\delta$ can be chosen as a projective subbundle such that 
$$
\delta \cdot \mathcal{X}_1 \; =\; \ell H,
$$
where $\ell H$ is viewed as a member of the complete linear system $|\mathcal{L}_1| = |L^{\otimes \ell }|$ 
on $\mathcal{X}_1 = X$.
For $\mathcal{X}$, we choose its proper $T$-equivariant desinguralization
$$
\iota : \hat{\mathcal{X}} \to \mathcal{X}
$$
 isomorphic over $\mathcal{X}\setminus \mathcal{X}_0$. 
Put $\hat{\pi} := \pi \circ \iota$.
Consider the $T$-invariant irreducible reduced divisor $\mathcal{H}$ on $\hat{\mathcal{X}}$ obtained 
as the closure in $\hat{\mathcal{X}}$ of the preimage of
$$
\bigcup_{t\in \Bbb C^*}\;\{t\}\times \psi (t) H
$$
under the mapping $\iota$, where $H$ on $X$ is viewd as a subset $\Bbb P^*(V_{\ell})$ via the Kodaira embedding
$ X \subset  \Bbb P^*(V_{\ell})$ associated to the complete linear system $|L^{\otimes \ell }|$.
Then we have the following equality of divisors on $\hat{\mathcal{X}}$:
$$
\iota^* (\delta \cdot \mathcal{X})\; = \;  \hat{D} + \ell \mathcal{H},
\leqno{(5.4)}
$$
where $\hat{D}$ is an effective divisor on $\mathcal{X}$ with support sitting in $\mathcal{X}_0$ set-theoretically.
Since $\mathcal{H}$ is a $T$-invariant divisor on $\hat{\mathcal{X}}$, the $T$-action on $\hat{\mathcal{X}}$ lifts 
to a $T$-linearization of $\hat{\mathcal{Q}} := \mathcal{O}_{\hat{\mathcal{X}}}(\mathcal{H})$.
Since $\mathcal{L} = \mathcal{O}_{\mathcal{X}}(\delta \cdot \mathcal{X})$, by (5.4), we obtain
$$
\hat{\mathcal{L}} \; =\; \mathcal{O}_{\hat{\mathcal{X}}}(\hat{D}) \otimes \hat{\mathcal{Q}}^{\otimes \ell}.
\leqno{(5.5)}
$$
For the direct image sheaf $F:= \hat{\pi}_*\hat{\mathcal{Q}}$ over $\Bbb A^1$, let $F_z$ be the fiber of $F$ 
over each $z \in \Bbb A^1$. Then we have a $T$-equivariant rational map
$$
\eta : \hat{\mathcal{X}} \to \Bbb P^* (F) 
$$
whose restriction over each $z \in \Bbb A^1\setminus \{0\}$ is the Kodaira embedding $\eta_z : \hat{\mathcal{X}}_z 
\hookrightarrow \Bbb P^*(F_z)$ associated to the complete linear system  $|\hat{\mathcal{Q}}_z|$ on $\hat{\mathcal{X}}_z$.
Put $\mathcal{Y}_z := \eta_z (\hat{\mathcal{X}}_z)$. Then the open subset $\hat{\pi}^{-1}(\Bbb A^1 \setminus \{0\})$ of 
$\hat{\mathcal{X}}$ is naturally identified with the $T$-invariant subset
$$ 
{\mathcal{Y}}^{\circ} :=\bigcup_{0\neq z\in \Bbb A^1}
{\mathcal{Y}}_z
$$ 
of $\Bbb P^*(F)$.
Let ${\mathcal{Y}}$ be the ${T}$-invariant 
subvariety of $\Bbb P^*(F)$ obtained 
as the closure  of 
${\mathcal{Y}}^{\circ}$ in $\Bbb P^*(F)$, i.e., ${\mathcal{Y}}$ is the meromorphic image of 
$\hat{\mathcal{X}}$ 
under the rational map  $\eta$.
Then the restriction 
$$
{\pi}^{}_{\mathcal{Y}} : \;{\mathcal{Y}}\;\to\; \Bbb A^1
$$  
to ${\mathcal{Y}}$ 
of the natural projection of $\Bbb P^*(F)$ onto $\Bbb A^1$ is a 
$T$-equivariant projective morphism 
with a relatively very ample invertible sheaf 
$$
\mathcal{Q}\; :=\;
 \mathcal{O}^{}_{\Bbb P^*(F)}(1)_{|{\mathcal{Y}}}
$$ 
on the fiber space ${\mathcal{Y}}$ over $\Bbb A^1$.
Note that $\hat{\pi} = {\pi}^{}_{\mathcal{Y}} \circ \eta$.
The ${T}$-action on $\hat{\mathcal{Q}}$ naturally induces a ${T}$-action on $F$, and 
it then induces a
${T}$-action on $\mathcal{O}_{{\mathcal{Y}}/\Bbb A^1}(1)$ 
covering the ${T}$-action on ${\mathcal{Y}}$.
By the affirmative solution of $T$-equivariant Serre's conjecture, we have a $T$-equivariant trivialization
$$
F \; \cong \; \Bbb A^1 \times F_0,
$$
where this isomorphism can be chosen in such a way that the Hermitian metric $\rho_1\, (=\,\rho_{\ell}\,{}_{|\ell = 1})$ as in the introduction
on  
$$
F_1\; = \;V_1\; =\; H^0(X,L)
$$ 
is taken to a Hermitian metric on $F_0$ which is preserved by the action of the compact subgroup $S^1 \subset T$ 
 (see \cite{D3}).
By this trivialization, $F_0$ can be identified with $F_1\,( = V_1)$, so that the $T$-action on $F_0$ induces
a representation
$$
\varphi : \; T \to \operatorname{GL}(V_1).
$$
Hence $({\mathcal{Y}}, \mathcal{Q}, \varphi )$ is a test configuration for $(X,L)$ of exponent 1.
Since $\hat{\mathcal{Q}} = \mathcal{O}_{\hat{\mathcal{X}}}(\mathcal{H})$,  the base point set $B$ for the subspace of 
$H^0(\hat{\mathcal{X}}_0, \hat{\mathcal{Q}}_0)$ 
associated to $F_0$ contains no components of dimension $n$.
However, replacing $\hat{\mathcal{X}}$ by its suitable birational model obtained from $\hat{\mathcal{X}}$ by a sequence of 
${T}$-equivariant blowing-ups with centers sitting over $B$, we may assume without loss of generality that $B$ is purely $n$-dimensional, i.e., $B = \emptyset$.
 Now the rational map
 $\eta : \hat{\mathcal{X}}\to {\mathcal{Y}} \subset \Bbb P^*(F)$
is holomorphic, and hence
$$
\hat{\mathcal{Q}}\; =  \; \eta^*\mathcal{Q},
$$ 
as required. This together with (5.5) completes the proof of Theorem 5.2.

\medskip\noindent
{\em Remark \/$5.6$.} Note that the divisor $\hat{D}$ on $\hat{\mathcal{X}}$ is preserved by the $T$-action.
Since $\mathcal{O}_{\hat{\mathcal{X}}}(\hat{D}) = 
\eta^*\mathcal{Q}^{\otimes \ell}\otimes \hat{\mathcal{L}}^{-1}$, the actions of  $\,T\, (= \Bbb C^*)$ on 
 $\mathcal{Q}$
and $\hat{\mathcal{L}}$ induce a $T$-action on 
the invertible sheaf $\mathcal{O}_{\hat{\mathcal{X}}}(\hat{D})$. 
Let $\zeta$ be a natural nonzero section for $\mathcal{O}_{\hat{\mathcal{X}}}(\hat{D})$ on $\hat{\mathcal{X}}$ 
having  $\hat{D}$ as the divisor $\operatorname{zero}(\zeta )$ of the zeroes.
Then the action of each element $t$ of $T$ on the line $\Bbb C \zeta$ is written as  
$$
\zeta \;\mapsto \; t^{\alpha} \zeta,
$$
where $\alpha \in \Bbb Z$ is the weight of the $T$-action on $\Bbb C \zeta$. 

\medskip
For test configurations $\mu$ and $\bar{\mu}:= (\mathcal{Y}, \mathcal{Q}, \varphi )$ 
above, each $t \in T$ 
not as a complex number but as an element 
of the  group $T$ of transformation on $\mu$ and $\bar{\mu}$
will be written as $g^{}_{\mu}(t)$ and $g^{}_{\bar{\mu}}(t)$, repectively. 
Let $\operatorname{Aut}(\hat{\mathcal{L}})$ and $\operatorname{Aut}(\mathcal{Q})$
denote the 
groups of all biholomorphisms of 
the total spaces of $\hat{\mathcal{L}}$ and $\mathcal{Q}$, respectively.
Then for $\varphi$ in Theorem 5.2, 
the $T$-linearization of $\mathcal{Q}$ defines a $T$-action on 
the real line bundle $|\mathcal{Q}|^2 := \mathcal{Q}
\otimes\bar{\mathcal{Q}}$ over $\mathcal{X}$
by
$$
g^{}_{\bar{\mu}}(t)\cdot |q |^2 := |g^{}_{\bar{\mu}}(t)\cdot q |^2 = |\tilde{\varphi}(t)(q)|^2,
\qquad (t, q ) \in T\times \mathcal{Q},
$$
where $\tilde{\varphi}: T \to  \operatorname{Aut}(\mathcal{Q})$
denotes the homomorphism induced by $\varphi$.
Note also that the $T$-linearization of $\hat{\mathcal{L}}$ induces a $T$-action on the real line bundle $|\hat{\mathcal{L}}|^2 := \hat{\mathcal{L}}
\otimes\bar{\mathcal{L}}$ such that
$$
g^{}_{\mu}(t)\cdot |\sigma |^2 := |g^{}_{\mu}(t)\cdot \sigma |^2 = |\tilde{\psi}(t)(\sigma )|^2,
\qquad (t, \sigma ) \in T\times \hat{\mathcal{L}},
$$
where $\tilde{\psi}: T \to \operatorname{Aut}(\hat{\mathcal{L}})$ denotes the homomorphism induced by 
$\psi$.  Note that both $g^{}_{\bar{\mu}}(t)$ and $g^{}_{\mu}(t)$ come from the same $T$-action.
Then for $\hat{\mathcal{Q}}:= \eta^*\mathcal{Q}$, by Theorem 5.2, we see that
$$
|\hat{\mathcal{L}}|^{2/\ell} \; =\; |\zeta |^{2/\ell}\, |\hat{\mathcal{Q}}|^2,
\leqno{(5.7)}
$$
where $T_{\Bbb R}$ acts on the real line $\Bbb R |\zeta |^{2/\ell}$ with weight $2\alpha/\ell$, so that
$g_{\mu} (t)\cdot |\zeta |^{2/\ell} = t^{2\alpha/\ell} |\zeta |^{2/\ell}$ for all $t \in T_{\Bbb R}$.
Since birational morphisms $\iota$ and $\eta$ are isomorphic over $\Bbb A^1\setminus\{ 0\}$,
by restricting them to $\{ z\neq 0\}$,
we can identify the line bundles $\hat{\mathcal{L}}$ and $\hat{\mathcal{Q}}$  
with $\mathcal{L}$ and $\mathcal{Q}$, respectively. Hence  (5.7) restricts to
$$
|\mathcal{L}|^{2/\ell} \; = \; |\zeta|^{2/\ell} |\mathcal{Q}|^2, 
\qquad z\neq 0.
\leqno{(5.8)}
$$
\smallskip\noindent
{\em Remark \/$5.9$.} 
The restriction of $\zeta$ to $z= 1$  gives a non-vanishing holomorphic section for 
$\mathcal{O}_{\hat{\mathcal{X}}}(\hat{D}){}_{|\hat{\mathcal{X}}_0}$. Define a Hermitian metric $\rho$
for $\mathcal{O}_{\hat{\mathcal{X}}}(\hat{D}){}_{|\hat{\mathcal{X}}_0}$ by
$$
|\zeta_{|\hat{\mathcal{X}}_0}|^2_{\rho} = 1
$$
everywhere on $\hat{\mathcal{X}}_0$. Then  by Theorem 5.2, when restricted to $z=1$, 
we may assume that $\mathcal{L}$ and $\mathcal{Q}^{\otimes \ell}$ coincides holomorphically and metrically.
In particular, any Hermitian metric for $L$ can be viewed as a Hermitian metric for 
$\mathcal{Q}_{|\hat{\mathcal{X}}_0}$ via the identification of $\hat{\mathcal{X}}_0$ with $X$.

\bigskip\noindent
{\footnotesize
{\sc Department of Mathematics}\newline
{\sc Osaka University} \newline
{\sc Toyonaka, Osaka, 560-0043}\newline
{\sc Japan}}
\end{document}